\documentclass[english,12pt]{amsart}
\usepackage[english]{babel}
\usepackage{amsmath,amssymb,enumerate,amsthm}
\usepackage[latin1]{inputenc}
\usepackage[T1]{fontenc}
\usepackage{cite}
\usepackage{latexsym}

\newtheorem{lemma}{Lemma}[section]
\newtheorem{teorema}[lemma]{Theorem}
\newtheorem{prop}[lemma]{Proposition}

\newtheorem{lema}[lemma]{Lemma}
\newtheorem{rk}[lemma]{Remark}

\newcommand{\PGL}{\mathrm{PGL}(2,\mathbb{Z})}

\newcommand{\R}{\mathbb{R}}
\newcommand{\Q}{\mathbb{Q}}
\newcommand{\Z}{\mathbb{Z}}
\newcommand{\rf}{\rfloor}
\newcommand{\lf}{\lfloor}

\def\={\;=\;}
\def\.={\;\dot{=}\;}
\usepackage{fancyhdr}

\begin{document}
\title{On a theorem of Serret on continued fractions}
\author{Paloma Bengoechea}
\address{Department of Mathematics, University of York, York, YO10 5DD, United Kingdom}
\email{paloma.bengoechea@york.ac.uk}

\begin{abstract} 
A classical theorem in continued fractions due to Serret shows that for any two irrational numbers $x$ and $y$ related by a transformation $\gamma$ in $\PGL$ there exist $s$ and $t$ for which the complete quotients $x_s$ and $y_t$ coincide. In this paper we give an upper bound in terms of $\gamma$ for the smallest indices $s$ and $t$.
\end{abstract}

\maketitle

A classical theorem of Serret \cite{S} states that the standard continued fraction expansions of two irrational numbers $x=[n_0,n_1,\ldots]$ and $y=[m_0,m_1,\ldots]$ are the same after a finite number of steps if and only if $x$ and $y$ are $\PGL$-equivalent. Serret's theorem does not quantify how large the smallest indices $s$ and $t$ have to be such that the complete quotients $x_s=[n_s,n_{s+1},\ldots]$ and $y_t=[m_t,m_{t+1},\ldots]$ coincide. The more complicated the matrix $\gamma$ relating $x$ and $y$ is, the larger $s$ or $t$ becomes. A quantification of this fact, namely bounding $s$ and $t$ in terms of $\gamma$, is the goal of this paper. The main theorem, in Section 2, gives an upper bound for the smallest $s$ and $t$ such that $x_s=y_t$. The bound we give is the best possible in the sense that for each $\gamma\in\PGL$ one can construct many real numbers $x$ for which $s$ or $t$ almost attain the bound, and for many matrices $\gamma$ the bound is actually attained.
In order to prove the theorem, in the first section  we introduce the set $\Gamma(x)$ of transformations of $\PGL$ related to the continued fraction expansion of the real number $x$. We characterize it by certain linear inequalities in Proposition 1.2.   Sets of transformations related to other continued fraction algorithms have been characterized in \cite{K} in a similar way. 


\section{Linear inequalities characterizing the continued fraction transformations}

We denote by $\Gamma$ the group $\PGL$ and by $\varepsilon$ and $T$ the transformations
$$
\varepsilon=\begin{pmatrix} 0 &1\\1 &0\end{pmatrix},\qquad
T=\begin{pmatrix} 1 &1\\ 0 &1\end{pmatrix}
$$
 that correspond to the inversion and the translation for the usual action on the projective line
$$
\begin{pmatrix}a &b\\c &d\end{pmatrix}x:=~\dfrac{ax+b}{cx+d}.
$$ 

The positive continued fraction of a real number $x$, 
$$
x=n_0+\dfrac{1}{n_1+\dfrac{1}{n_2+\dfrac{1}{\ddots}}}\qquad (n_0\in\mathbb{Z},\, n_i\geq 1\ \forall i\geq 1),
$$
also denoted by $x=[n_0, n_1,\ldots]$, is given by the algorithm
\begin{align*}
&x_0=x,\\
&n_i=\lf x_i\rf,\quad  x_{i+1}=\dfrac{1}{x_i-n_i}\=\varepsilon T^{-n_i}(x_i)\quad\mbox{if $x_i\not\in\Z$}\qquad(i\geq 0).
\end{align*} 
If $x_i\in\Z$, the algorithm stops after $n_i=\lfloor x_i\rfloor$.
Clearly each complete quotient $x_i=[n_i,n_{i+1},\ldots]$ is the image of $x$ by a matrix $\gamma_i=\gamma_{i,x}\in\Gamma$, given explicitly by 
$$
\gamma_0:=\mathrm{Id},\qquad \gamma_i:=\begin{pmatrix} 0 &1\\1 &-n_{i-1}\end{pmatrix}\cdots\begin{pmatrix} 0 &1\\1 &-n_0\end{pmatrix}\qquad(i\geq 1)
$$
and recursively by
\begin{equation}\label{matrices}
\gamma_{0}=\mathrm{Id},\qquad \gamma_{i+1}=\varepsilon T^{-n_i}\gamma_{i}\qquad(i\geq 0).
\end{equation}
A key idea is to replace the sequence $(\gamma_1,\gamma_2,\gamma_3,\ldots)$ of elements in $\Gamma$ by the unordered set 
$$
\Gamma(x)=\left\{\gamma_1,\gamma_2,\gamma_3,\ldots\right\}\subset\Gamma.
$$
When $x\in\R\backslash\Q$, the set $\Gamma(x)$ is infinite, whereas it is finite if $x\in\Q$. If $x=[n_0,\ldots,n_t]$, then 
$$
\Gamma(x)=\left\{\gamma_1,\ldots,\gamma_{t+1}\right\}.
$$

The $i$-th convergent of $x$ is denoted by $\dfrac{p_i}{q_i}=[n_0,\ldots,n_i]$.
The integers $p_i$ and $q_i$ satisfy the recurrence 
$$
p_{-2}=0, \qquad p_{-1}=1,\qquad p_i=n_ip_{i-1}+p_{i-2} \qquad (i\geq 0),
$$
$$
q_{-2}=1,\qquad q_{-1}=0, \qquad q_i=n_iq_{i-1}+q_{i-2} \qquad (i\geq 0),
$$
the equation
$$
p_{i+1}q_i-p_iq_{i+1}=(-1)^i
$$
and the inequalities
\begin{equation}\label{ineq 1}
q_i\geq q_{i-1}\geq0 \quad (i\geq 0), \qquad\quad q_i>q_{i-1}>0\quad (i\geq 2),
\end{equation}
\begin{equation}\label{ineq 2}
|p_i|\geq|p_{i-1}| \quad(i\geq2), \qquad\quad |p_i|>|p_{i-1}|\quad (i\geq4).
\end{equation}

It is well known that any rational number $p/q$ satisfying
$$
\left|\dfrac{p}{q}-x\right|<\dfrac{1}{2q^2}
$$ 
is a convergent of $x$.

The numbers $\delta_i$ ($i\geq -1$) defined by 
$$
\delta_i=(-1)^i(p_{i-1}-q_{i-1}x)
$$
satisfy the recurrence
$$
\delta_{-1}=x,\qquad \delta_{0}=1,\qquad
\delta_{i+1}=-n_i\delta_{i}+\delta_{i-1},
$$
and the inequalities $1=\delta_{0}>\delta_1>\ldots\geq0$. If $x$ is rational, then $x_i=p_i/q_i$ for some $i$ and the recurrence stops with $\delta_{i+1}=0$; if $x$ is irrational, the $\delta_i$ are all positive and converge to 0 with exponential rapidity. 
One has
$$
n_i=\left\lfloor\dfrac{\delta_{i-1}}{\delta_{i}}\right\rfloor \quad\mbox{if $\delta_i\neq 0$},
$$
and
$$
\gamma_i^{-1}=\begin{pmatrix}p_{i-1} &p_{i-2}\\q_{i-1} &q_{i-2}\end{pmatrix},\qquad \gamma_i\begin{pmatrix}x\\ 1\end{pmatrix}=\begin{pmatrix}\delta_{i-1}\\ \delta_i\end{pmatrix}.
$$

\begin{lema}\label{lema convergentes} If $r/t$ and $s/u$ are two rational numbers such that $t,u>~0$,\,  $\dfrac{r}{t}\leq x\leq\dfrac{s}{u}$ and $ru-st=\pm1$, then $r/t$ or $s/u$ is a convergent of $x$.
\end{lema}

\textbf{Proof.} Let $x$ be an element in $\Big[\dfrac{r}{t},\dfrac{s}{u}\Big]$. If one of the inequalities 
$$
\Big|\dfrac{r}{t}-x\Big|<\dfrac{1}{2t^2} \quad\mbox{or}\quad  \Big|\dfrac{s}{u}-x\Big|<\dfrac{1}{2u^2}
$$
is satisfied, then $r/t$ or $s/u$ is a convergent of $x$. If not, then
$$
\dfrac{1}{tu}=\Big|\dfrac{r}{t}-\dfrac{s}{u}\Big|=\Big|\dfrac{r}{t}-x\Big|+\Big|\dfrac{s}{u}-x\Big|\geq\dfrac{1}{2t^2}+\dfrac{1}{2u^2}.
$$
The inequality above can only hold if $t=u=1$, in which case $s=r+1$, so $\Big[\dfrac{r}{t},\dfrac{s}{u}\Big]=[r,r+1]$, and then either $\dfrac{s}{u}=x=\dfrac{p_0}{q_0}$ or $\dfrac{r}{t}=\lf x\rf=\dfrac{p_0}{q_0}$.
\begin{flushright}
$\square$
\end{flushright} 

Throughout, for any pair of real numbers $a,b$, we denote by $|a,b|$ the interval $[a,b]$ if $a\leq b$ or $[b,a]$ if $b<a$.


\begin{prop} \label{prop 1}
For any real number $x$, the set $\Gamma(x)$ equals $W-(W_1\cup W_2)$, where
$$
W=\left\{\gamma\in\Gamma\, \mid -1\leq\gamma(\infty)\leq 0,\, \gamma(x)>1\right\}
$$
and
$$
W_1=\left\{\gamma\in W\, \mid\, \gamma(\infty)=0,\, \det(\gamma)=1\right\},
$$ 
$$
W_2=\left\{\gamma\in W\, \mid\, \gamma(\infty)=-1,\, \mathrm{det}(\gamma)=-1\right\}.
$$
The sets $W_1$ and $W_2$ have respectively exactly one and at most one element. (The inequality $\gamma(x)>1$ in the definition of $W$ includes $\gamma(x)=\infty$.)
\end{prop}

\textbf{Proof.} One easily checks that
\begin{align*}
&W_1\=\left\{\begin{pmatrix}0 &-1\\1 &-1-n_0\end{pmatrix}\right\},
\\
&W_2\=\left\{\begin{array}{cl}\left\{
\begin{pmatrix}-1 &1+n_0\\1 &-n_0\end{pmatrix}\right\} &\mbox{if $x\in\Z$ or $n_1\geq2$},\\
\emptyset &\mbox{if $n_1=1$}.
\end{array}\right.
\end{align*}
This proves the second statement. It is also easy to see that $\Gamma(x)\subseteq W-(W_1\cup W_2)$. Indeed, $\gamma_{i}(x)=\dfrac{\delta_{i-1}}{\delta_{i}}>1$ and $\gamma_{i}(\infty)=-\dfrac{q_{i-2}}{q_{i-1}}$ for $i\geq 1$ together with the inequalities \eqref{ineq 1} imply $\gamma_i\in W-(W_1\cup W_2)$.

Therefore we only have to show that $W\subseteq\Gamma(x)\cup W_1\cup W_2$.
Let $\gamma\in\Gamma$ satisfy 
$$-1\leq\gamma(\infty)<0,\qquad\gamma(x)>1.$$ 
Let $$\gamma^{-1}=\begin{pmatrix}r &s\\t &u\end{pmatrix},\qquad\gamma=\begin{pmatrix}u &-s\\-t &r\end{pmatrix}.$$
The conditions $\gamma(\infty)<0$ and $\gamma(x)>0$ imply that $u$ and $t$ have the same sign, as well as $ux-s$ and $-tx+r$, so $x\in\Big|\dfrac{r}{t},\dfrac{s}{u}\Big|$. By Lemma \ref{lema convergentes}, $r/t$ or $s/u$ is a convergent of $x$. 

If $r/t$ is a convergent of $x$, then $\gamma$ and $\gamma^{-1}$ are of the form
\begin{equation}\label{caso 1}
\gamma^{-1}=\begin{pmatrix} p_i &p_{i-1}+kp_i\\q_i &q_{i-1}+kq_i\end{pmatrix},\qquad\gamma=\begin{pmatrix} q_{i-1}+kq_i &-p_{i-1}-kp_i\\-q_i &p_i\end{pmatrix}
\end{equation}
or
\begin{equation}\label{caso 2}
\gamma^{-1}=\begin{pmatrix} p_i &-p_{i-1}-kp_i\\q_i &-q_{i-1}-kq_i \end{pmatrix},\qquad\gamma=\begin{pmatrix} -q_{i-1}-kq_i &p_{i-1}+kp_i\\-q_i &p_i \end{pmatrix}
\end{equation}
with $k\in\mathbb{Z}$ and $i\geq 0$. In the case \eqref{caso 1}, since $q_{i-1}\leq q_i$, we have that $-1\leq\gamma(\infty)=-\dfrac{q_{i-1}}{q_i}-k<0$ if and only if $k=0$ and $i\geq 1$, or $k=1$ and $i=0$. Moreover, if $x\not\in\Z$, then $\gamma(x)=\dfrac{\delta_{i}}{\delta_{i+1}}-k>1$ if and only if $k\leq n_{i+1}-1$ because $\left\lfloor\dfrac{\delta_{i}}{\delta_{i+1}}\right\rfloor=n_{i+1}$. Hence $-1\leq\gamma(\infty)<0$ and $\gamma(x)>1$ if and only if $k=0$ and $i\geq 1$ or $x\in\Z$, $k=1$ and $i=0$, or $n_1\geq 2$, $k=1$ and $i=0$. In the first case
$$
\gamma=\begin{pmatrix}
q_{i-1} &-p_{i-1}\\
-q_i &p_i
\end{pmatrix}=\gamma_{i+1}\qquad (i\geq 1);
$$
in the other two cases
$$
\gamma=\begin{pmatrix}-1 &1+n_0\\1 &-n_0\end{pmatrix}\in W_2.
$$

In the case \eqref{caso 2}, if $-1\leq\gamma(\infty)=\dfrac{q_{i-1}}{q_i}+k<0$, then $k<0$ and $\gamma(x)=-\dfrac{\delta_{i}}{\delta_{i+1}}+k<~0$, so $\gamma\not\in W$.
\\ 

If $s/u$ is a convergent of $x$, then $\gamma$ and $\gamma^{-1}$ are of the form
\begin{equation}\label{caso 3}
\gamma^{-1}=\begin{pmatrix}p_{i-1}+kp_i &p_i\\q_{i-1}+kq_i &q_i\end{pmatrix},\qquad\gamma=\begin{pmatrix}q_i &-p_i\\-q_{i-1}-kq_i &p_{i-1}+kp_i\end{pmatrix}
\end{equation}
or
\begin{equation}\label{caso 4}
\gamma^{-1}=\begin{pmatrix}-p_{i-1}-kp_i &p_i\\-q_{i-1}-kq_i &q_i\end{pmatrix},\qquad \gamma=\begin{pmatrix}q_i &-p_i\\q_{i-1}+kq_i &-p_{i-1}-kp_i\end{pmatrix}                                                            
\end{equation}
with $k\in\mathbb{Z}$ and $i\geq 0$. In the case \eqref{caso 3}, $\gamma(x)=\dfrac{\delta_{i+1}}{\delta_{i}-k\delta_{i+1}}>1$ if and only if $k=n_{i+1}$. If $k=n_{i+1}$, then $-1\leq\gamma(\infty)=-\dfrac{q_i}{q_{i+1}}<0$. Hence
$$
\gamma=\begin{pmatrix} q_i &-p_i\\ -q_{i+1} &p_{i+1}\end{pmatrix}=\gamma_{i+2}\qquad (i\geq 0).
$$

In the case \eqref{caso 4}, if $-1\leq\gamma(\infty)=\dfrac{q_i}{q_{i-1}+kq_i}<0$, then $k<0$ and $\gamma(x)=\dfrac{\delta_{i+1}}{-\delta_{i}+k\delta_{i+1}}\leq0$, so $\gamma\not\in W$.
\\

Finally, let $\gamma=\begin{pmatrix}0 &1\\1 &u\end{pmatrix}$ be the unique element of $W$ satisfying $\gamma(\infty)=0$ and $\mathrm{det}(\gamma)=~-1$. We have $1<\gamma(x)=\dfrac{1}{x+u}$ if and only if $0<x+u<1$, namely $u=-\lfloor x\rfloor$, and $\gamma=\gamma_{1}\in\Gamma(x)$.
\begin{flushright}
$\square$
\end{flushright}


\section{Main theorem}

We recall the classic theorem of Serret about continued fraction expansions of equivalent numbers (see \cite{S}):
\begin{teorema}[Serret]  \label{Hurwitz}
Two irrational numbers $x$ and $y$ are $\Gamma$-equivalent if and only if there exist $s,t\geq 0$ such that $x_s=y_t$.
\end{teorema}

Theorem \ref{Hurwitz} (or its proof) does not give any bound independent of $x$ and $y$ for the indexes $s$ and $t$. The result we state below gives the best possible bound for the index $s$ in terms of the matrix relating $x$ and $y$ (see Remark \ref{rk}). 

We introduce the convention that $\infty$ has zero partial quotients in the rest of the paper.

\begin{teorema}\label{teo} 
Let $\gamma\in\Gamma$ and let $r$ be the number of partial quotients of $\gamma^{-1}(\infty)$. For every real number $x$ there exists an index $s\leq r+3$ such that the complete quotient $x_s$ is also a complete quotient of $\gamma(x)$.
\end{teorema}

Of course the bound for $t$ is obtained in the same way in terms of the number of partial quotients of $\gamma(\infty)$.
\\

\textbf{Proof.} 
Let $y=\gamma(x)$.
Since 
$$
x_s=\gamma_{s,x}(x)\qquad\mbox{and}\qquad y_t=\gamma_{t,y}(y)=(\gamma_{t,y}\gamma)(x),
$$
we must prove that there exists $s\leq r+3$ such that $\gamma_{s,x}=\gamma_{t,y}\gamma$ for some $t\geq 1$. In other words, we prove next that, 
if $\gamma_{i,x}\neq\gamma_{t,y}\gamma$, then $i\leq r+2$. Let $i\geq 0$ be such that $\gamma_{i,x}\gamma^{-1}\not\in\Gamma(y)$.
Proposition \ref{prop 1} implies one of the following conditions:
\begin{equation}\label{conditions}
\gamma_{i,x}\gamma^{-1}(\infty)\geq0\qquad\mbox{or}\qquad
\gamma_{i,x}\gamma^{-1}(\infty)\leq-1.
\end{equation}
Suppose $\gamma^{-1}(\infty)\neq\infty$ and set $p/q=\gamma^{-1}(\infty)$ with $(p,q)=1$.
If the first inequality in \eqref{conditions} is satisfied, we have
$$
\gamma_{i,x}\Big(\dfrac{p}{q}\Big)\=\dfrac{q_{i-2}p-p_{i-2}q}{-q_{i-1}p+p_{i-1}q}\, \geq\, 0,
$$ 
which implies $\dfrac{p}{q}\in\left|\dfrac{p_{i-2}}{q_{i-2}},\dfrac{p_{i-1}}{q_{i-1}}\right|$. Then by Lemma \ref{lema convergentes} we have that either $p_{i-1}/q_{i-1}$ or $p_{i-2}/q_{i-2}$ is a convergent of $p/q$.  
Similarly, if the second inequality in \eqref{conditions} is satisfied, we have 
$$
\gamma_{i,x}\Big(\dfrac{p}{q}\Big)\=\dfrac{q_{i-2}p-p_{i-2}q}{-q_{i-1}p+p_{i-1}q}\, \leq\, -1,
$$ 
which implies $\dfrac{p}{q}\in\left|\dfrac{p_{i-1}-p_{i-2}}{q_{i-1}-q_{i-2}},\dfrac{p_{i-1}}{q_{i-1}}\right|$. It is easy to see that $\dfrac{p_{i-1}-p_{i-2}}{q_{i-1}-q_{i-2}}\in\left|\dfrac{p_{i-3}}{q_{i-3}},\dfrac{p_{i-1}}{q_{i-1}}\right|$. Hence by Lemma \ref{lema convergentes} we have that either $p_{i-1}/q_{i-1}$ or $p_{i-3}/q_{i-3}$ is a convergent of $p/q$. 
Thus, if we denote by $r$ the number of partial quotients (or convergents) of $p/q$, we have that $i-3\leq r-1$.

Now suppose $\gamma^{-1}(\infty)=\infty$. Then we have that, for all $i\geq 1$,
$$
\gamma_{i,x}\gamma^{-1}(y)=\gamma_{i,x}(x)\qquad\mbox{and}\qquad \gamma_{i,x}\gamma^{-1}(\infty)=\gamma_{i,x}(\infty),
$$
so 
$$
\gamma_{i,x}\gamma^{-1}(y)>1\qquad\mbox{and}\qquad  -1\leq\gamma_{i,x}\gamma^{-1}(\infty)\leq 0.
$$
If, moreover, $\gamma_{i,x}\gamma^{-1}(\infty)\not\in\left\{-1,0\right\}$, then, by Proposition \ref{prop 1}, we have that $\gamma_{i,x}\gamma^{-1}\in\Gamma(y)$. Now, if $\gamma_{i,x}(\infty)=0$, then $\gamma_{i,x}=\gamma_{1,x}=\begin{pmatrix} 0 &1\\1 &-n_0\end{pmatrix}$. If $\gamma_{i,x}(\infty)=~-1$, then $n_1=1$ and $\gamma_{i,x}=\gamma_{2,x}=\begin{pmatrix} 1 &-n_0\\-1 &1+n_0\end{pmatrix}$. Therefore $i\leq 2$.
\begin{flushright}
$\square$
\end{flushright}

\begin{rk}\label{rk} Let $\gamma\in\Gamma$ and set $\gamma^{-1}(\infty)=[n_0,\ldots,n_{r-1}]$.
One can construct many real numbers $x$ for which $x_{r+1}$ is not a complete quotient of $\gamma(x)$ by letting $x=[n_0,\ldots,n_{r-2},n_{r-1}-1,1,n_{r+1},\ldots]$ where $n_{r+i}$ are arbitrary for all $i\geq 1$. Indeed, for such an $x$, we have that $\gamma^{-1}(\infty)=p_r/q_r$, so
$\gamma_{r+1,x}\gamma^{-1}(\infty)=\infty$, and then $\gamma_{r+1,x}\not\in\Gamma(y)\gamma$. Therefore $x_{r+1}$ is not a complete quotient of $\gamma(x)$. Moreover, $\gamma_{r+2,x}\gamma^{-1}(\infty)=0$, and if $\mathrm{det}(\gamma)=\mathrm{det}(\gamma_{r+2,x})$, then  $\gamma_{r+2,x}\not\in\Gamma(y)\gamma$ and $x_{r+2}$ is not a complete quotient of $\gamma(x)$ either.
\end{rk}



\end{document}